\documentclass[12pt]{amsart}

\usepackage{fullpage}
\usepackage{amssymb}
\usepackage{amsmath}
\usepackage[colorlinks, breaklinks, linkcolor=blue]{hyperref}
\usepackage{caption}

\newcommand{\F}{\mathbb F}
\newcommand{\Z}{\mathbb Z}

\def\b{\mathfrak b}
\def\u{\mathfrak u}

\newtheorem*{thm*}{Theorem}

\title
{The modality of a Borel subgroup in a simple algebraic group of type $E_8$}

\author[S.~M.~Goodwin]{Simon M. Goodwin}
\address{School of Mathematics, University of Birmingham, Birmingham, B15
2TT, United Kingdom}  \email{s.m.goodwin@bham.ac.uk}

\author[P.~Mosch]{Peter Mosch}
\author[G.\ R\"ohrle]{Gerhard R\"ohrle}
\address{Fakult\"at f\"ur Mathematik, Ruhr-Universit\"at Bochum, D-44780 Bochum, Germany}
\email{peter.mosch@rub.de} \email{gerhard.roehrle@rub.de}


\begin{document}

\begin{abstract}
Let $G$ be a simple algebraic group over
an algebraically closed field $k$, where
$\mathrm{char}\, k$ is either 0 or a good prime for $G$.
We consider the {\em modality} $\mathrm{mod}(B : \u)$ of
the action of a Borel subgroup $B$ of $G$ on the Lie algebra $\u$ of the
unipotent radical of $B$,
and report on computer calculations used to show that $\mathrm{mod}(B:\u) = 20$,
when $G$ is of type $\mathrm E_8$.
This completes the determination of the values for 
$\mathrm{mod}(B:\u)$ for $G$ of exceptional type.
\end{abstract}

\maketitle

\setcounter{section}{1}

Let $G$ be a simple algebraic group over
an algebraically closed field $k$, where
$\mathrm{char}\, k$ is either 0 or a good prime for $G$.
Let $B$ be a Borel subgroup of $G$ with unipotent radical $U$,
and let $\u = \mathrm{Lie}\, U$.  The modality of the adjoint action of the
Borel subgroup $B$ on $\u$ is defined by
$
\mathrm{mod}(B : \u) := \max_{i \in \Z_{\ge 0}} ( \dim(\u(B,i)) - i),
$
where $\u(B,i) := \{ x \in \u \mid \dim B \cdot x = i\}$, and
is intuitively the maximal number of parameters on which a family of $B$-orbits in $\u$ depends.
As $\mathrm{mod}(B:\u)$ is an important invariant of the action of $B$ on $\u$ it
is of interest to determine its value.
As a general reference for the modality of the action of an algebraic group
on an algebraic variety, we refer to \cite{Vi}.

In \cite{JR} the values of $\mathrm{mod}(B:\u)$ are determined for $G$ up to rank 7 excluding
type $\mathrm B_7$ and $\mathrm E_7$; they are found by combining lower bounds
from \cite[Prop.\ 3.3]{Roe} with upper bounds obtained by computer calculation, and are presented
in \cite[Tables II and III]{JR}. We refer to the introduction of \cite{JR}
and the references therein for prior history of finding values of $\mathrm{mod}(B:\u)$, and for
motivation. The known
values of $\mathrm{mod}(B:\u)$ were extended to $G$ up to rank 8 excluding type $\mathrm E_8$
in \cite[\S5]{GMR2}; this
required computer calculations, which are explained below.  We note also that,
as is explained in \cite[\S6]{GG}, results in \cite{PS} can be used to determine
$\mathrm{mod}(B:\u)$ for $G$ of type $\mathrm A_l$ for $l \le 15$.

Our main theorem gives $\mathrm{mod}(B:\u)$ in the case $G$ is of type $\mathrm E_8$.

\begin{thm*}
Let $G$ be a simple algebraic group of type $\mathrm E_8$ over
the algebraically closed field $k$, where
$\mathrm{char}\, k = 0$ or $\mathrm{char}\, k > 5$.
Let $B$ be a Borel subgroup of $G$ with unipotent radical $U$.
Then $\mathrm{mod}(B : \u) = 20$.
\end{thm*}

Our theorem
completes
the list of values for $\mathrm{mod}(B : \u)$ for $G$ of exceptional type as presented in the table below.

\begin{table}[ht!b]
\renewcommand{\arraystretch}{1.5}
\begin{tabular}{r|cccccc}\hline
Type of $G$ & $G_2$ & $F_4$ & $E_6$ & $E_7$ & $E_8$ \\
\hline
$\mathrm{mod}(B : \u)$ & 1 & 4 & 5 & 10 & 20  \\
\hline
\end{tabular}
\caption*{Modality of the action of $B$ on $\u$}
\end{table}

We move on to review the computer programme from \cite{GMR2} and explain how it
was adapted to show that 20 is an upper bound for $\mathrm{mod}(B:\u)$ for $G$ of type $\mathrm E_8$.
Thanks to \cite[Prop.\ 3.3]{Roe} it is known that 20 is a lower bound for $\mathrm{mod}(B:\u)$,
so combining these bounds proves our theorem.  In the discussion below we refer to
$\mathrm{mod}(U:\u)$ and $\mathrm{mod}(U:\u^*)$, which are defined analogously to $\mathrm{mod}(B:\u)$.

It is shown in \cite[Prop.\ 5.4]{Go} that each $U$-orbit in $\u$ admits a so called
\emph{minimal representative}.  As explained in \cite[\S2]{GR1}, the minimal
representatives are partitioned into certain locally closed subsets $X_c$ of $\u$ for $c$
running over some index set $C$.  This gives a parametrization of the
$U$-orbits in $\u$, so we can deduce that $\mathrm{mod}(U:\u) = \max_{c \in C} \dim X_c$,
and thus by \cite[Thm.\ 5.1]{GMR2} that $\mathrm{mod}(B:\u) = \max_{c \in C} \dim X_c - \mathrm{rank}\, G$.
An algorithm for determining all the varieties $X_c$ for $c \in C$ is given in \cite[\S3]{GR1}.
This algorithm was programmed in GAP, \cite{GAP}, and subsequent developments were made in \cite{GMR1} including
calls to SINGULAR, \cite{SIN}.  The resulting programme was used to obtain the parametrization
of the $U$-orbits in $\u$ when $G$ is of rank up to 7 except for type $\mathrm E_7$; so this
can also be
used to determine $\mathrm{mod}(B:\u)$ in these cases.

The results in \cite[\S3]{GMR2} show that a similar algorithm is valid
for the coadjoint action of $U$ on $\u^*$.  In particular, there is a
parametrization of minimal representatives of $U$-orbits in $\u^*$ by
certain locally closed subsets $\mathcal X_c$ of $\u^*$ for $c$
running over an index set $C$.  This algorithm was programmed
and used to obtain a complete description of the varieties $\mathcal X_c$, when $G$ has rank up to 8, with the
exception of type $\mathrm E_8$. Since $\mathrm{mod}(U:\u) = \mathrm{mod}(U:\u^*)$, see \cite[Thm~1.4]{Roe2},
we have $\mathrm{mod}(B:\u) = \max_{c \in C} \dim \mathcal X_c - \mathrm{rank}\, G$.  Thus this allowed
us to determine $\mathrm{mod}(B:\u)$ when $G$ has rank up to 8, with the
exception of type $\mathrm E_8$.

The algorithm for determining the varieties $\mathcal X_c$ for $c$ in $C$ involves a certain
polynomial-resolving subroutine, as explained in \cite[\S3]{GMR1}.
This is the most complicated and computationally expensive part of the programme.
We adapted our algorithm, so that in cases where the programme is not
able to resolve all the polynomial conditions in a specified amount of time it
simply disregards these unresolved conditions.  Thus the modified computation determines a variety
$\mathcal Y_c \supseteq \mathcal X_c$, which we can view as an ``upper bound''
for a parametrization of the minimal representatives in $\mathcal X_c$,
so that $\dim \mathcal Y_c \ge \dim \mathcal X_c$,  for each $c \in C$.
Consequently, $\mathrm{mod}(U:\u) \le \max_{c \in C} \dim \mathcal Y_c$.

We ran the programme for the case $G$ of type $\mathrm E_8$ and determined a variety
$\mathcal Y_c$ for every $c$ in $C$.  From the output of the computation we obtain that
$\mathrm{mod}(B:\u) \le 20$ as required to verify our theorem.

\smallskip

We move on to mention consequences of our calculations for the finite groups
of rational points, when $G$ is defined over a finite field.
Suppose that $G$ is defined and split over the field $\F_p$ where
$p$ is a good prime $p$ for $G$.
Let $q$ be a power of $p$ and denote by $G(q)$ the group of $\F_q$-rational points of $G$.
Also assume that $B$ is defined over $\F_q$, so
$U$ is defined over $\F_q$ and $U(q)$ is a Sylow $p$-subgroup of $G(q)$.
We write $k(U(q))$ for the number of conjugacy classes of $U(q)$ (which is also
the number of complex irreducible characters of $U(q)$).
As explained in \cite[\S4]{GMR2}, the parametrization of the coadjoint
$U$-orbits in $\u^*$ by the varieties $\mathcal X_c$ for $c$ in $C$, gives
a method to calculate $k(U(q))$.  In fact in the cases considered in \cite{GMR2},
there is a polynomial $g(t) \in \Z[t]$ such that $k(U(q)) = g(q)$; and, moreover,
$g(t)$ does not depend on $p$.  Our adapted programme
 calculates a polynomial $h(t) \in \Z[t]$ such that $k(U(q)) \le h(q)$ and $h(t)$ does not
 depend on $p$.
Moreover, an upper bound for $\mathrm{mod}(U:\u)$ can be easily read off
as the degree of $h(t)$; we refer to \cite[\S5]{GMR2} for further details.
Note that we do not claim here that $k(U(q))$ is necessarily a polynomial in $q$
for $G$ of type $\mathrm E_8$, and remark
that \cite[Thm.\ 1.4]{PS} suggests that this might not be the case for general $G$.

\smallskip

We end by noting that our calculation of $\mathrm{mod}(B:\u)$ can be used to
determine the dimension of the commuting varieties of $\u$ and $\b$ as are studied in
\cite{GR2} and \cite{GG}, respectively.

\end{document}